\documentclass[a4paper]{amsart}

\usepackage{amsfonts}
\usepackage{amsmath}
\usepackage{amssymb}
\usepackage{amsthm,color}
\usepackage{enumerate}
\usepackage{mathtools}
\mathtoolsset{showonlyrefs}
 \usepackage{transparent}
 \usepackage{stackengine,tikz}

\usepackage{color,graphicx,enumerate,wrapfig,amssymb,mathtools}
\usepackage{epsfig}
\usepackage{pxfonts}
\usepackage{graphicx}

\usepackage{eucal}

\newcommand{\R}{{\mathbb R}}

\newcommand{\e}{\varepsilon}
\newcommand{\al}{\alpha}
\newcommand{\be}{\beta}
\newcommand{\de}{\delta}
\newcommand{\ka}{\kappa}
\newcommand{\la}{\lambda}
\newcommand{\si}{\sigma}

\newcommand{\dist}{\operatorname{dist}}

\newcommand{\loc}{\operatorname{loc}}

\newcommand{\D}{\nabla}
\newcommand{\p}{\partial}

\usepackage{amsmath}
\usepackage{amsfonts}
\usepackage{caption}
\usepackage{subcaption}
\usepackage{esint}

\newtheorem{theorem}{Theorem}
\theoremstyle{plain}

\newtheorem{remark}{Remark}

\numberwithin{equation}{section}

\makeatletter
\@namedef{subjclassname@2020}{%
  \textup{2020} Mathematics Subject Classification}
\makeatother

\begin{document}

\author{Layan El Hajj}
\address[Layan El Hajj]
{Mathematics Division \newline 
\indent  American University in Dubai \newline 
\indent Dubai, United Arab Emirates} 
\email[Layan El Hajj]{lhajj@aud.edu}

\author{Seongmin Jeon}
\address[Seongmin Jeon]
{Department of Mathematics \newline 
\indent KTH Royal Institute of Technology \newline 
\indent 100 44 Stockholm, Sweden} 
\email[Seongmin Jeon]{seongmin@kth.se}

\author{Henrik Shahgholian}
\address[Henrik Shahgholian]
{Department of Mathematics \newline 
\indent KTH Royal Institute of Technology \newline 
\indent 100 44 Stockholm, Sweden} 
\email[Henrik Shahgholian]{henriksh@kth.se}

\title[ Symmetry in free boundary problems]
{Symmetry for a fully nonlinear free boundary  problem with highly singular term }

\date{\today}
\thanks{HS is supported by Swedish Research Council.}

\begin{abstract}
    In this paper we prove radial symmetry for solutions to a
     free boundary problem with a singular right hand side, in both elliptic and parabolic regime. 
     More exactly, in the unit ball 
     $B_1$ we consider a solution to the fully nonlinear elliptic problem
    $$
     \begin{cases}
    F(D^2u)=f(u)&\text{in }B_1 \cap \{u >0 \},\\
    u=M&\text{on }\p B_1,\\
    0\le u<M&\text{in }B_1,
\end{cases}
$$
where the right hand side $f(u) $, near $u=0$,  behaves like $u^a$ with negative values for $a \in (-1,0)$.
Due to lack of $C^2$-smoothness of  both $u$ and  the free boundary $\partial\{u>0\}$, we cannot apply the well-known  Serrin-type boundary point lemma. We circumvent this by an exact  assumption on a first order 
 expansion  and the decay on the second order, along with  an ad-hoc comparison principle.

We treat equally  the parabolic case of the problem, and state a corresponding result.

\end{abstract}

\maketitle 

\tableofcontents

\section{Introduction}

\subsection{Background }

 In this paper we deal with rigidity question of  a fully nonlinear Dirichlet   problem in the unit ball with constant boundary values. Our  problem,     due to   presence of a discontinuous  right hand side   gives rise to a free boundary; the r.h.s. in our problem is given by $f(u) \approx  u^a $ with $ -1 < a < 0$,  when $0 < u \approx 0$.
Such singularities    imply much  less regularity of the solutions up to the free boundary, as is a customary assumptions in rigidity theory. A further difficulty in our problem is the fact that the  moving plane technique, which is  the main tool, usually requires Lipschitz right hand sides. 
This is naturally missing in our problem, due to the highly singular r.h.s. 
   
   The fully nonlinear case, as dealt with here, introduces yet another technical problem, which has to do with the equation itself and what definition should be  appropriate. An obvious choice of definition would be viscosity (which we also adopt in the bulk) but the problem is the regularity issue for such solutions at the free boundary.
   
   In the Laplacian case,    problems with highly singular r.h.s., and in a general setting,
    have  been studied earlier by Alt-Phillips \cite{AltPhi86}, as minimizers,  and by  Teixeira \cite{Tei18} as viscosity solutions.     See also a recent work by Desilva-Savin \cite{DeSSav21} for the case $-3< a <-1$, where the authors consider minimizers of a corresponding functional and  define viscosity solutions through an   asymptotic behaviour for solutions  near the free boundary.
    
Two of the current authors \cite{ElHSha21} have also treated the symmetry problem for the Laplacian case, where the authors considered (a stronger form of) asymptotic expansion, for the solution, to by pass the low-regular nature of the solutions and that of the free boundary. 

The fully nonlinear case with singular terms has been considered earlier; see e.g., \cite{AraTei13}, \cite{FelQuaSir12}.
However, no results has been found in the literature where symmetry aspects have been considered, for singular right hand side.

 It would be an interesting project to consider analysis of the free boundary problem in the nonlinear setting, in the way Desilva-Savin \cite{DeSSav21} have done for the Laplacian case. Also the weaker version of asymptotic expansion of DeSilva-Savin is far from sufficient for treating any symmetry problems.
    To circumvent this obstacle,   we shall assume a second order  asymptotic expansion, which is justified\footnote{Although this is a challenging problem  even in the Laplacian case, it is reasonable to expect such an expansion.} from a regularity point of view; see Asymptotic Property \eqref{eq:asymp-ell}, \eqref{eq:asymp-par-x}.

  We will now define our problem more specifically, in both elliptic and parabolic cases.
  The elliptic case of our problem is given by the equation (with $\partial \{u > 0\}$ a priori unknown)
 \begin{align}\label{eq:sol-PDE-FN}
\begin{cases}
    F(D^2u)=f(u)&\text{in }B_1 \cap \{u >0 \},\\
    u=M&\text{on }\p B_1,\\
    0\le u<M&\text{in }B_1,
\end{cases}
\end{align}
where  $B_1 \subset \R^n$ ($n \geq 2$) is the unit ball,    $M>0$ is a given fixed constant, and the functions $f$ and $F$ are defined below.

The parabolic counter part of our problem is expressed through the equation
\begin{align}\label{eq:sol-PDE-FN-par}
\begin{cases}
    F(D^2u)-\p_tu=f(u)&\text{in }Q_1  \cap \{u >0 \} ,\\
    u(x,t)=M_t&\text{on } \p_pQ_1,\\
    0\le u(x,t)<M_t&\text{in }Q_1.
\end{cases}
\end{align}
Here, $Q_1=B_1\times(0,1]$ is the parabolic unit cylinder and $\p_pQ_1=(\p B_1\times[0,1])\cup(B_1\times\{0\})$ is its parabolic boundary. $M_t$ is a positive continuous function of $t\in[0,1]$.

In both problems we shall have an additional assumption of the behaviour of the solution close to the (unknown)
 free boundary $\partial \{u > 0\}$. This assumption is expressed in terms of an asymptotic expansion, that assumes the solution close to $C^1$ free boundaries behave like a one-dimensional solution and the second term in the expansion has much faster decay with respect to the distance to the free boundary. 
 
For the Laplacian case Alt-Phillips considered existence and regularity of solutions and the free boundary in their work \cite{AltPhi86}, through minimization, which does not work for non-divergence equations. 
 For the fully nonlinear case one may consider a singular perturbation technique, as in \cite{AraTei13}. The latter approach seems plausible for the parabolic case.


\subsection{Main Results}
To specify the function $f$ in \eqref{eq:sol-PDE-FN} and \eqref{eq:sol-PDE-FN-par}, we fix $-1<a<0$, $\ka_0>0$ and $\e_1>0$ small. We assume that $f:\R\to\R$ satisfies the following properties
\begin{align}
    \label{eq:assump-f-FN}
    \begin{cases}
        (a)\, \sup_{\e>0}\frac{f(\rho+\e)-f(\rho)}\e\le\ka_0\rho^{a-1},\quad \rho>0,\\
        (b)\,f(\rho)=0,\quad \rho\le0, \\
        (c)\, f(\rho)=\rho^a,\quad 0<\rho<\e_1,\\
        (d)\,|f(\rho)|\le C_0,\quad \rho\ge \e_1.
    \end{cases}
\end{align}
To describe the operator $F$, let $S=S(n)$ be the space of $n\times n$ symmetric matrices, $\Lambda\ge\la>0$ be constants, and $P^-_{\la,\Lambda}$, $P^+_{\la,\Lambda}$ be the extremal Pucci operators defined by

\begin{equation*}
P^-_{\la,\Lambda}(N)=\la\Sigma_{e_i>0}e_i+\Lambda\Sigma_{e_i<0}e_i,\quad P^+_{\la,\Lambda}(N)=\la\Sigma_{e_i<0}e_i+\Lambda\Sigma_{e_i>0}e_i,
\end{equation*}
where each $e_i$ is an  eigenvalue of $N\in S$. Our assumption on $F:S\to\R$ is the following:
\begin{align}\label{cond:F}
    \begin{cases}
     \text{- } F\text{ is convex}.\\
    \text{- } F(0)=0.\\
    \text{- } F\text{ is homogeneous of degree }1,\text{ i.e., }F(rN)=rF(N)\text{ for any }r>0,\,N\in S. \\
    \text{- }\text{There are constants }\Lambda\ge\la>0\text{ such that for any } N_1, N_2\in S, \hbox{we have}\\ \qquad P_{\la,\Lambda}^-(N_1-N_2)\le F(N_1)-F(N_2)\le P_{\la,\Lambda}^+(N_1-N_2).\\
    \text{- }1\le\Lambda/\la\le\eta_0\text{ for some }\eta_0>1,\text{ defined in Appendix~\ref{appen:comp-prin}.}\\
    \text{- $F$ is Hessian, i.e., $F(R^TNR)=F(N)$ for each orthogonal matrix $R$ and $N\in S$.}
    \end{cases}
\end{align}

\medskip

\noindent \underline{Asymptotic Property} 

For $-1<a<0$ we fix the following constant throughout the paper: 
\begin{align*}
\be=\frac2{1-a}.
\end{align*}

\begin{itemize}

\item[]
\underline{Elliptic case:} \quad
We say that a function $u:B_1\to\R$ satisfies \emph{Asymptotic Property} if at any free boundary point $x^0\in\p\{u>0\}$ the  asymptotic expansion holds \begin{align}\label{eq:asymp-ell}
|u(x)-A_{x^0}((x-x^0)\cdot\nu^{x^0})_+^\be|\le C_0|x-x^0|^{2+\de_\be},\quad x\in B_{d_{x_0}}(x^0).
\end{align} 
Here, $C_0>0$ and $0<\de_\be<\be-1$ are fixed constants independent of $x^0$, while a unit vector 
$\nu^{x^0}\in\p B_1$,  a radius $d_{x_0}\ge c_0>0$ and a constant $A_{x^0}>0$ depend on $x^0$.\\

\item[]
\underline{Parabolic case:}
\quad 
The parabolic case needs a similar type of asymptotic property. Let $u:Q_1\to \R$ be a function, with $0<T_1<1$ the first time we encounter a free boundary point $\p\{u>0\}$. We say that $u$ satisfies \emph{(Parabolic) Asymptotic Property} if at every free boundary point $z^0=(x^0,t^0)\in\p\{u>0\}$ with $t^0>T_1$ the asymptotic expansion holds
\begin{equation}
    \label{eq:asymp-par-x}
    \left|u(x,t)-A_{z^0}((x-x^0)\cdot\mu^{z^0})^\be_+\right|\le C_0\left(|x-x^0|+\sqrt{|t-t^0|}\right)^{2+\de_\be},\quad (x,t)\in\tilde Q_{d_{z^0}}(z^0).
    \end{equation}
Here, $C_0>0$ and $0<\de_\be<\be-1$ are fixed constants independent of $z^0$, while a spatial unit vector $\mu^{z^0}\in\p B_1$ in $\R^n$, a constant $A_{z^0}>0$ and a radius 
$$d_{z^0}:=\min\{|1-t^0|^{1/2},\,|t^0-T_1|^{1/2}  \}$$
depend on $z^0$.

\end{itemize}

\begin{remark}\label{rem:FB-smooth} 
It is noteworthy that the above asymptotic property \eqref{eq:asymp-ell} implies that the free boundary is uniformly $C^{1,|a|}$ in the elliptic case. 

On the other hand, in the parabolic case the asymptotic property \eqref{eq:asymp-par-x} implies 
 $C^{1/2, |a|/2}$ in  time direction and $C^{1,|a|}$ in space directions for each point on $\p\{u>0\}\cap\{t>T_1\}$, with uniform norm, depending on the distance of the free boundary point to the time slice $\{t=T_1\}$.

\end{remark}

The asymptotic expansions \eqref{eq:asymp-ell} and \eqref{eq:asymp-par-x} guarantee that the homogeneous rescalings $u_r(x):=\frac{u(rx+x^0)}{r^\be}$ converge  (up to a subsequence) to $q_{x^0}(x):=A_{x^0}(x\cdot\nu^{x^0})^\be_+$, and similarly $u_r(x,t):=\frac{u(rx+x^0,r^2t+t^0)}{r^\be}$ converge to (a time-independent version) $q_{z^0}(x,t):=A_{z^0}(x\cdot\mu^{z^0})_+^\be$. Since $f(u)=u^a$ near the free boundary $\p\{u>0\}$, it is reasonable to expect that they solve $F(D^2q_{x^0})=q_{x^0}^a$ and $F(D^2q_{z^0})-\p_tq_{z^0}=q_{z^0}^a$, which will determine $A_{x^0}$ and $A_{z^0}$, respectively. See Appendix~\ref{sec:appen-A} for their exact values and properties.

By the result of Remark~\ref{rem:FB-smooth} on the regularity of the free boundary $\p\{u>0\}$, it is easy to see that the vector $\nu^{x^0}$ in \eqref{eq:asymp-ell} should be the unit normal to the free boundary at $x^0$, which points toward $\{u>0\}$. Similarly, in \eqref{eq:asymp-par-x}, if $\nu^{z^0}=(\nu^{z^0}_x,\nu^{z^0}_t)$ is the unit normal to $\p\{u>0\}$ pointing toward $\{u>0\}$ then $\mu^{z^0}=\frac{\nu^{z^0}_x}{|\nu^{z^0}_x|}$.

Our results are the following theorems.

\begin{theorem}\label{thm:sym-FN}
For $-1<a<0$, let $u$ be a viscosity solution of \eqref{eq:sol-PDE-FN}.
Suppose that $f$ satisfies \eqref{eq:assump-f-FN} and $f(\rho)\ge0$ for $\rho\ge M-\e_1$, and that the asymptotic property \eqref{eq:asymp-ell} holds. Then $u$ is spherically symmetric in $B_1$ and $\p_{|x|}u\ge 0$.
\end{theorem}

Our second main result concerning the parabolic equation is as follows.

\begin{theorem}\label{thm:sym-FN-par}
For $-1<a<0$, let $u$ be a viscosity solution of \eqref{eq:sol-PDE-FN-par}. Suppose that $f$ satisfies \eqref{eq:assump-f-FN}, $f(\rho)\ge0$ for $\rho\ge \inf_{0\le t\le1}M_t-\e_1$ and $\inf_{\e>0}\frac{f(\rho+\e)-f(\rho)}\e\ge-\ka_0\rho^{a-1}$ for $\rho>0$. We assume that the asymptotic property \eqref{eq:asymp-par-x} holds. Then, for each $0<t<1$, $u(\cdot,t)$ is spherically symmetric around the origin, and 
$\p_{|x|}u (x,t)\ge 0$.
\end{theorem}

The assumption (in parabolic case) that close to every free boundary point we have AP property is somehow restrictive, and we believe the theorem  should be true without this. The trouble point (which we could not deal with) is that if there are more points failing \eqref{eq:asymp-par-x} after time $T_1$, then standard comparison/maximum principles are much harder to work out at such points, as their appears several possibilities, some of which we could handle but some of them resisted all our attempts, and would need some further tools, unknown to us at this stage.

We remark that when $a\ge 0$ the above symmetry results easily follow from the standard moving plane technique, due to the monotonicity of the r.h.s. $f(u)=u^a$ close to $u=0$. 
We have thus left out this particular case, even though a similar machinery of moving plane would be needed, but the subtle technical approach in this paper is not needed anymore for this case.


\subsection{Notation}
For $z^0=(x^0,t^0)\in\R^{n+1}$, where $x^0\in\R^n$ and $t\in\R$, and $r>0$, we let
 \begin{align*}
    B_r(x^0)&=\{x\in\R^n\,:\,|x-x^0|<r\}\,:\,\text{Euclidean ball}, \\
    Q_r(z^0)&=B_r(x^0)\times(t^0,t^0+r^2]\,:\,\text{parabolic cylinder}, \\
    \tilde Q_r(z^0)&=B_r(x^0)\times(t^0-r^2,t^0+r^2)\,:\,\text{full parabolic cylinder}, \\
    \p_pQ_r(z^0)&=\overline{Q_r(z^0)}\setminus Q_r(z^0)\,:\,\text{parabolic boundary}.
\end{align*}
For a general bounded open set $\Omega\subset\R^{n+1}$, its parabolic boundary $\p_p\Omega$ is defined as the set of all points $z\in\p \Omega$ such that for any $\e>0$ the interior of the parabolic cylinder $Q_\e(z)$ contains at least one  point in 
$ \Omega^c$.


\section{Proof of Theorem~\ref{thm:sym-FN} }

For $0\le \tau<1$ and $x=(x_1,\cdots,x_n)$, we define \begin{align*}
    \begin{cases}
        \Pi_\tau=\{x\in\R^n\,:\,x_1=\tau\}\,:\,\text{the hyperplane},\\
        x^\tau=(2\tau-x_1,x_2,x_3,\cdots\,x_n)\,:\,\text{the reflection of }x\text{ with respect to }\Pi_\tau,\\
        \Sigma_\tau=\{x\,:\,x^\tau\in B_1,\,x_1<\tau\}\,:\,\text{the reflection of the right side of $B_1$},\\
        u_\tau(x)=u(x^\tau)\,:\,\text{the reflection of }u\text{ with respect to }\Pi_\tau.
     \end{cases}   
\end{align*}
By Hessian property of $F$ it is enough to prove symmetry only in $x_1$-direction. It  further suffices 
  to show that $u_\tau\ge u$ in $\Sigma_\tau$ for any $0<\tau<1$. The theorem will then follow from this.

Before proving $u_\tau\ge u$ in $\Sigma_\tau$, $0<\tau<1$, we observe that (by Hessian property)
 $u_\tau$ satisfies the equation \begin{align}
    \label{eq:u_tau-pde}
    F(D^2u_\tau)=f(u_\tau)\quad\text{in }\Sigma_\tau\cap\{u_\tau>0\},
\end{align}
and also  the asymptotic property at every $y\in\p\{u_\tau>0\}$ \begin{align}
    \label{eq:u_tau-asymp}
    |u_\tau(x)-A_{y^\tau}\left((x-y)\cdot\tilde\nu^{y^\tau}\right)_+^\beta|\le C_0|x-y|^{2+\de_\be},\quad x\in B_{d_{y^\tau}}(y).
\end{align}
Here, $\tilde\nu^{y^\tau}$ is the reflection of the unit normal $\nu^{y^\tau}$ of $u$ at $y^\tau$ through the plane $\{x_1=0\}$ (i.e., if $\nu^{y^\tau}=(\nu_1,\cdots,\nu_n)$ then $\tilde\nu^{y^\tau}=(-\nu_1,\nu_2,\cdots,\nu_n)$). Notice that $\tilde\nu^{y^\tau}$ is equal to the unit normal of $u_\tau$ at $y$ pointing into $\{u_\tau>0\}$. Indeed, for \eqref{eq:u_tau-pde}, one can compute, using $u_\tau(x)=u(x^\tau)=u(2\tau-x_1,x_2,\cdots,x_n)$, $D^2u_\tau(x)=R^TD^2u(x^\tau)R$ for a diagonal matrix $R$ with entries $-1,1,\cdots,1$. Then
$$
F(D^2u_\tau(x))=F(R^TD^2u(x^\tau)R)=F(D^2u(x^\tau))=f(u(x^\tau))=f(u_\tau(x))\quad\text{in }\Sigma_\tau\cap\{u_\tau>0\}.
$$
To prove \eqref{eq:u_tau-asymp}, we employ the asymptotic property of $u$ at $y^\tau\in\p\{u>0\}$ to deduce that for $x\in B_{d_{y^\tau}}(y)$ \begin{align*}
\left|u_\tau(x)-A_{y^\tau}\left((x-y)\cdot\tilde\nu^{y^\tau}\right)_+^\be\right|&=\left|u(x^\tau)-A_{y^\tau}\left((x^\tau-y^\tau)\cdot\nu^{y^\tau}\right)_+^\be\right|\\
&\le C_0|x^\tau-y^\tau|^{2+\de_\be}=C_0|x-y|^{2+\de_\be}.
\end{align*}

Using $P^+_{\la,\Lambda}(D^2u)\ge F(D^2u)-F(0)=f(u)\ge 0$ near $\p B_1$, we have by Hopf's Lemma (Proposition 2.6 in \cite{DaLSir07}) that \begin{align}\label{eq:hopf-bdry}
    \p_\nu u < 0\quad\text{on }\p B_1,
\end{align}
for any unit vector $\nu$ on $\p B_1$ pointing into $B_1$. Thus $\p_{x_1}u>0$ in a small neighborhood of $e_1=(1,0,0,\cdots,0)$. This gives that for any $\tau\in(0,1)$ close to $1$, $u_\tau>u$ in $\Sigma_\tau$. Hence we can start our moving plane and we define $\tau_0$ as the smallest value such that $u_\tau\ge u$ in $\Sigma_\tau$ for all $\tau_0<\tau<1$.

\indent Towards a contradiction, we assume that $\tau_0>0$, and take sequences $\e_i\searrow0$, $\tau_i=\tau_0-\e_i>0$ satisfying $D_i:=\{w_{\tau_i}<0\}\cap\Sigma_{\tau_i}\neq\emptyset$, where $w_{\tau_i}:=u_{\tau_i}-u$. For such $D_i$, we define $$D_0:=\{x^0\in\overline{B_1}\,:\,x^{i_k}\to x^0\text{ for some subsequence }x^{i_k}\in D_{i_k}\}\subset\overline{\Sigma_{\tau_0}}.$$
From $D_i\neq\emptyset$, we have $D_0\neq\emptyset$. For $w_{\tau_0}:=\lim_{i\to\infty}w_{\tau_i}=u_{\tau_0}-u$, we claim \begin{align}
    \label{eq:w-D_0}
    w_{\tau_0}=|\D w_{\tau_0}|=0\quad\text{in }D_0.
\end{align}
Indeed, the fact that $w_{\tau_i}<0$ in $D_i$ yields $w_{\tau_0}\le0$ on $D_0$, while the definition of $\tau_0$ implies $w_{\tau_0}\ge0$ in $\overline{\Sigma_{\tau_0}}\supset D_0$. Thus $w_{\tau_0}=0$ on $D_0$. To see that $\D w_{\tau_0}=0$ on $D_0$, let $x^0\in D_0$ and $x^i\in D_i$ with $x^i\to x^0$ over a subsequence. We also let $e\in\p B_1$ be an arbitrary direction. Since $D_i$ is open, there is a line segment $(z_1^i,z_2^i)\subset D_i$ such that $z_1^i$, $z_2^i\in\p D_i$, $z_1^i=z_2^i+re$ for some $r>0$, and $x^i$ lies on $(z^i_1,z^i_2)$. Since $w_{\tau_i}\ge 0$ on $\p \Sigma_{\tau_i}$, we have $w_{\tau_i}=0$ on $\p D_i$. In particular, $w_{\tau_i}(z_1^i)=w_{\tau_i}(z_2^i)=0$. Combining this with the $C^1$-regularity of $w_{\tau_i}$ in $\overline{\Sigma_{\tau_i}}$, we deduce that if $z_3^i\in(z_1^i,z_2^i)$ is a local minimum point of $w_{\tau_i}$ on $(z_1^i,z_2^i)$ then $\p_ew_{\tau_i}(z_3^i)=0$. Over a subsequence $[x^i,z_3^i]$ converges either to a point $\{x^0\}$ or to a line segment $[x^0,z^0]\subset D_0$. In either case we have $\p_ew_{\tau_0}(x^0)=0$. Since $e$ and $x^0$ are arbitrary, we conclude that $\D w_{\tau_0}=0$ on $D_0$.

\indent Now, we decompose $D_0$ into two parts $$
D_0=(D_0\cap\{u>0\})\cup(D_0\cap\{u=0\}),
$$
and will prove that these sets are empty, contradicting that $D_0\neq\emptyset$.\\

\medskip

\noindent
\underline{\emph{Claim  A:}}\quad   $D_0\cap\{u>0\}=\emptyset.$ \\

\emph{Claim A} will follow once we show the following: \\
(A1) $D_0\cap\Sigma_{\tau_0}\cap\{u>0\}=\emptyset$,\\
(A2) $D_0\cap \Pi_{\tau_0}\cap\{u>0\}=\emptyset.$\\

We first prove (A1). From  \eqref{eq:u_tau-pde} and (a) in \eqref{eq:assump-f-FN}, we can see that $w_{\tau_0}$ satisfies in $\Sigma_{\tau_0}\cap\{u>0\}$
 \begin{align}\label{eq:w-subsol-1}
P_{\la,\Lambda}^-(D^2w_{\tau_0})\le F(D^2u_{\tau_0})-F(D^2u)=f(u_{\tau_0})-f(u)\le \ka_0u^{a-1}w_{\tau_0}.
\end{align}
We also note that $w_{\tau_0}\ge 0$ in $\overline{\Sigma_{\tau_0}}$ and $u^{a-1}$ is bounded in every compact subset of $\Sigma_{\tau_0}\cap\{u>0\}$. By the strong minimum principle (Proposition~2.6 in \cite{DaLSir07}), for every connected component $C$ of $\Sigma_{\tau_0}\cap\{u>0\}$, either $w_{\tau_0}=0$ in $C$ or $w_{\tau_0}>0$ in $C$. We claim that \begin{align}
    \label{eq:w>0}
    w_{\tau_0}>0\quad\text{in }\Sigma_{\tau_0}\cap\{u>0\},
\end{align}
which readily implies (A1). To this end, we assume to the contrary that for some component $C$, 
$w_{\tau_0}=0$ in $C$. 
Since $C$ is an open set and 
$\p C\cap(\p\Sigma_{\tau_0}\setminus \Pi_{\tau_0})=\emptyset$, we can take two points $y,z$ with 
$y\in\p C\cap\Sigma_{\tau_0}$ and $z\in C$ satisfying $u(y)=u_{\tau_0}(y)=0$, $u(z)=u_{\tau_0}(z)>0$ and $z=y+re_1$ for some $r>0$. Then the reflected points $y^{\tau_0}$, $z^{\tau_0}\in B_1\cap \{x_1>\tau_0\}$ satisfy $u(y^{\tau_0})=0$, $u(z^{\tau_0})>0$ and $y=z+re_1$. This is a contradiction, since $u$ is nondecreasing in $x_1$-direction in $B_1\cap\{x_1>\tau_0\}$ by the definition of $\tau_0$; 
 we could also  have  used  the fact that   $u_{\tau}\ge u$ for $\tau=\frac{y^{\tau_0}_1+z^{\tau_0}_1}2>\tau_0$.\\

\indent To prove (A2), suppose there is a point $x^0\in D_0\cap \Pi_{\tau_0}\cap\{u>0\}$.
If $x^0\in \p B_1$, then $\p_{x_1}u(x^0)=-1/2\p_{x_1}w_{\tau_0}(x^0)=0$, which contradicts \eqref{eq:hopf-bdry}. Thus we may assume that $x^0\in B_1$. Then we can take a small ball in $\Sigma_{\tau_0}\cap\{u>0\}$ which touches $\Pi_{\tau_0}$ at $x^0$. From \eqref{eq:w>0}, $w_{\tau_0}(x^0)=0$ and $w_{\tau_0}>0$ in that ball. By using \eqref{eq:w-subsol-1} and Hopf's lemma, we get $$
\p_{x_1}w_{\tau_0}(x^0)<0.
$$
This is a contradiction since $\D w_{\tau_0}=0$ on $D_0$, and the proof of \emph{Claim A} is completed.\\


\medskip

\noindent
\underline{\emph{Claim  B:}} \quad  $D_0\cap\{u=0\}=\emptyset$. \\

For each $D_i=\{w_{\tau_i}<0\}\cap \Sigma_{\tau_i}$ defined in the beginning of the proof, we take a minimum point $x^i$ of $w_{\tau_i}$ in $D_i$ (or, equivalently, in $\Sigma_{\tau_i}$), so that $\D w_{\tau_i}(x^i)=0$. This is possible since $w_{\tau_i}<0$ in $D_i$ and $w_{\tau_i}\ge0$ on $\p D_i$. Then, over a subsequence $x^i\to x^0\in D_0$. By the result of \emph{Claim A}, $x^0\in D_0\cap\{u=0\}$. Note that we have from $w_{\tau_0}=u_{\tau_0}-u=M-u>0$ on $\p\Sigma_{\tau_0}\setminus\Pi_{\tau_0}$ and $w_{\tau_0}=0$ on $D_0$ that $D_0\subset \Sigma_{\tau_0}\cup \Pi_{\tau_0}$. We then have the following three possibilities: 
\begin{enumerate}[(i)]
\item $x^0\in D_0\cap \Sigma_{\tau_0}\cap\{u=0\}$.
\item $x^0\in D_0\cap\Pi_{\tau_0}\cap\{u=0\}$ and $\Pi_{\tau_0}$ is orthogonal to $\p\{u>0\}$ at $x^0$.
\item $x^0\in D_0\cap\Pi_{\tau_0}\cap\{u=0\}$ and $\Pi_{\tau_0}$ is non-orthogonal\footnote{The reader may ask why this case is handled in our problem, despite it being ignored in most (if not all) other similar problems. This depends on the use of  the   maximum principle in small domains, that  does not work (at least we do not know how) in our case due to the singularity in the r.h.s.}
 to $\p\{u>0\}$ at $x^0$.
\end{enumerate}

Before discussing the above three cases, we prove that $x^0\in\p\{u>0\}\cap\p\{u_{\tau_0}>0\}$.

Indeed, it follows from $u(x^i)>u_{\tau_i}(x^i)\ge 0$ and $u(x^0)=0$ that $x^0\in\p\{u>0\}$. In addition, from $u>0$ on $\p B_1$ we see $D_0\cap\{u=0\}\subset \Sigma_{\tau_0}\cup(\Pi_{\tau_0}\setminus\p B_1)$. This implies $x^i\in \Sigma_{\tau_0}$ for large $i$. Using $u_{\tau_0}\ge u$ in $\Sigma_{\tau_0}$, we have $u_{\tau_0}(x^i)\ge u(x^i)>0$. Combining those with $u_{\tau_0}(x^0)=w_{\tau_0}(x^0)+u(x^0)=0$ yields $x^0\in\p\{u_{\tau_0}>0\}$.

\medskip

We simultaneously consider the cases (i) and (ii), for they both follow from the application of the result in Appendix~\ref{appen:comp-prin}. We claim that the asymptotic property gives \begin{align}
    \label{eq:w-small}
    w_{\tau_0}(x)=O\left(|x-x^0|^{2+\de_\be}\right).
\end{align}
Indeed, \eqref{eq:w-small} is trivial for the case (ii) due to its orthogonality assumption. For the case (i), we use the asymptotic expansions \eqref{eq:asymp-ell} of $u$ at $x^0$ and \eqref{eq:u_tau-asymp} of $u_{\tau_0}$ at $x^0$, respectively, to have \begin{align*}
    u(x)&=A_{x^0}\left((x-x^0)\cdot\nu^{x^0}\right)_+^\be+O\left(|x-x^0|^{2+\de_\be}\right),\\
    u_{\tau_0}(x)&=A_{(x^0)^{\tau_0}}\left((x-x^0)\cdot\tilde\nu^{(x^0)^{\tau_0}}\right)_+^\be+O\left(|x-x^0|^{2+\de_\be}\right).
\end{align*}
The fact that $u_{\tau_0}\ge u$ in $\Sigma_{\tau_0}$ and $x^0\in \Sigma_{\tau_0}$ implies $u_{\tau_0}\ge u$ in a neighborhood of $x^0$, which yields $\nu^{x^0}=\tilde\nu^{(x^0)^{\tau_0}}$. This in turn implies, employing the equality $A_y=\left[\be(\be-1)F\left(\nu^y\otimes\nu^y\right)\right]^{-\be/2}$, $y\in\p\{u>0\}$, in Appendix~\ref{sec:appen-A}, that $A_{x^0}=A_{(x^0)^{\tau_0}}$. Thus, \eqref{eq:w-small} holds in Case (i) as well.

Next, we observe that for small $r>0$ \begin{align*}
    P^-_{\la,\Lambda}(D^2w_{\tau_0})&\le f(u_{\tau_0})-f(u)=u_{\tau_0}^a-u^a\le 0\quad\text{in }A_r:=\Sigma_{\tau_0}\cap\{u>0\}\cap B_{r}(x^0).
\end{align*}
To use the result in Appendix~\ref{subsec:comparison-ell}, let the $n$-dimensional cone $K^\e$ and the function $H^\e$ be as in Appendix~\ref{subsec:comparison-ell}, with small $\e>0$ satisfying \eqref{eq:comparison-L}. Thanks to the $C^{1,|a|}$-regularity of $\p\{u>0\}$, we can take possible rotations and translations on $K^\e$ and $H^\e$ to obtain a cone $K$ and a function $H$ satisfying for small $r>0$ \begin{align}
        &K\cap B_r(x^0)\subset A_r,\\
        &P^-_{\la,\Lambda}(D^2H)=0\quad\text{in }K\cap B_r(x^0), \quad H=0\quad\text{on }\p K\cap B_r(x^0),\\
        &\label{eq:H-nondeg-seq}
H(x^j)\ge|x^j-x^0|^{2+\de_\be/2}\,\text{ for a sequence }x^j\in K\cap B_r(x^0)\text{ with }|x^j-x^0|\searrow0.
\end{align}
Since $K\cap\p B_r(x^0)$ is compactly supported in $\{w_{\tau_0}>0\}$, the minimum of  $w_{\tau_0}$ on $K\cap\p B_r(x^0)$ is strictly positive. Thus we have for some constant $c_*>0$ $$
c_*H\le w_{\tau_0}\quad\text{on } K\cap\p B_r(x^0).
$$
This, together with the fact that $c_*H=0\le w_{\tau_0}$ on $\p K\cap B_r(x^0)$, gives $$
w_{\tau_0}-c_*H\ge0 \quad\text{on }\p (K\cap B_r(x^0)).
$$
Moreover, $$
P^-_{\la,\Lambda}(D^2(w_{\tau_0}-c_*H))\le P^-_{\la,\Lambda}(D^2w_{\tau_0})-c_*P^-_{\la,\Lambda}(D^2H)\le 0\quad\text{in }K\cap B_r(x^0).
$$
Therefore, we have by the maximum principle $w_{\tau_0}\ge c_*H$ in $K\cap B_r(x^0)$. This contradicts
 \eqref{eq:w-small}--\eqref{eq:H-nondeg-seq}.\\

\medskip

\indent Now we consider  case (iii). Let $\nu^0=(\nu_1^0,\cdots,\nu_n^0)$ be the unit normal vector of $\p\{u>0\}$ at $x^0$ pointing toward $\{u>0\}$, i.e., $\nu^0=\nu^{x_0}$. Since the non-orthogonality assumption in (iii) gives $\nu_1^0\neq0$ and the fact that $w_{\tau_0}\ge0$ in $\Sigma_{\tau_0}$ implies $\nu_1^0\ge0$, we have $\nu_1^0>0$.

For each $i$ we take a free boundary point $y^i\in\p\{u>0\}$ closest to $x^i$. Then, for $\rho_i:=|x^i-y^i|$, $\nu^i:=\frac{x^i-y^i}{\rho_i}$ is the unit normal to $\p\{u>0\}$ at $y^i$. Due to the $C^1$-smoothness of $\p\{u>0\}$, $\nu^i\to\nu^0$, and thus $\nu^i_1\ge \nu^0_1/2>0$ for large $i$. Similarly, for each $i$ we take a point $\bar y^i\in\p\{u_{\tau_i}>0\}$ closest to $x^i$. Then, for $\bar\rho_i:=|x^i-\bar y^i|$, $\bar\nu^i:=\frac{x^i-\bar y^i}{\bar\rho_i}$ is the unit normal to $\p\{u_{\tau_i}>0\}$ at $\bar y^i$. $\nu^0_1>0$ implies that the first components $y^i_1$, $\bar y^i_1$ of $y^i$, $\bar y^i$, respectively, satisfy $y^i_1<\tau_i<\bar y^i_1$. This, combined with $x^i\in \Sigma_{\tau_i}\subset\{x_1<\tau_i\}$, yields $\bar\rho_i\ge\rho_i$ and $\bar\nu^i_1\le 0$. Notice that $x^i$, $y^i$, $\bar y^i$ and $(\bar y^i)^{\tau_i}$ all converge to $x^0$.

We claim that $\frac{\bar\rho_i}{\rho_i}\le 2$ for large $i$. Otherwise, $\frac{\bar\rho_{i}}{\rho_{i}}>2$ over a subsequence. Using $\rho_i=\dist(x^i,\p\{u>0\})$ and applying asymptotic property \eqref{eq:asymp-ell} of $u$ at $y^i$ gives for large $i$\begin{align*}
    u(x^i)\le A_{x^i}\rho_i^\be+C_0\rho_i^{2+\de_\be}\le \sqrt{3/2}A_{x^0}\rho_i^\be+C_0\rho_i^{2+\de_\be}\le\sqrt 2A_{x^0}\rho_i^\be.
\end{align*}
Similarly, using $\bar\rho_i=\dist(x^i,\p\{u_{\tau_i}>0\})$ and applying asymptotic property \eqref{eq:u_tau-asymp} of $u_{\tau_i}$ at $\bar y^i$, we obtain for large $i$ \begin{align*}
    u_{\tau_i}(x^i)&\ge A_{(\bar y^i)^{\tau_i}}\left((x^i-\bar y^i)\cdot\tilde\nu^{(\bar y^i)^{\tau_i}}\right)_+^\be-C_0|x^i-\bar y^i|^{2+\de_\be}=A_{(\bar y^i)^{\tau_i}}\bar\rho_i^\be-C_0\bar\rho_i^{2+\de_\be}\ge \frac{A_{x_0}}{\sqrt{2}}\bar\rho_i^\be.
\end{align*}
Here, in the second step we have used the fact that $\tilde\nu^{(\bar y^i)^{\tau_i}}$ is the unit normal $\bar\nu^i$ of $u_{\tau_i}$ at $\bar y^i$. If follows that $$
w_{\tau_{i}}(x^{i})=u_{\tau_{i}}(x^{i})-u(x^{i})\ge \frac{A_{x^0}}{\sqrt{2}}\bar\rho_{i}^\be-\sqrt{2}A_{x^0}\rho_{i}^\be>0,
$$
which contradicts $x^{i}\in D_{i}=\{w_{\tau_{i}}<0\}\cap\Sigma_{\tau_{i}}$.

Next, we consider the rescalings $$
q_i(x):=\frac{u(\rho_ix+y^i)}{\rho_i^\be},\qquad \bar q_i(x):=\frac{u_{\tau_i}(\rho_ix+y^i)}{\rho_i^\be}.
$$
Using the asymptotic property of $u$ at $y^i$ and $u_{\tau_i}$ at $\bar y^i$ again, \begin{align*}
    &|q_i(x)-A_{y^i}(x\cdot\nu^i)_+^\be|\le C_0\rho_i^{2+\de_\be-\be}|x|^{2+\de_\be},\\
    &\left|\bar q_i(x)-A_{(\bar y^i)^{\tau_i}}\left(\left(x+\frac{y^i-\bar y^i}{\rho_i}\right)\cdot\bar\nu^i\right)_+^\be\right|\le C_0\rho_i^{2+\de_\be-\be}\left|x+\frac{y^i-\bar y^i}{\rho_i}\right|^{2+\de_\be}.
\end{align*}
For large $i$ these estimates hold in $B_{1/2}(\nu^i)$, and thus in $B_{1/4}(\nu^0)$ as well. From $\frac{\bar\rho_i}{\rho_i}\le 2$, we see that $$
\frac{|y^i-\bar y^i|}{\rho_i}\le\frac{|y^i-x^i|+|x^i-\bar y^i|}{\rho_i}=\frac{\rho_i+\bar\rho_i}{\rho_i}\le 3.
$$
Thus, over a subsequence, \begin{align*}
    &q_i\to q_0(x):=A_{x^0}(x\cdot\nu^0)_+^\be\quad\text{uniformly in }B_{1/4}(\nu^0),\\
    &\bar q_i\to\bar q_0(x):=A_{x^0}((x+z^0)\cdot\bar\nu^0)_+^\be\quad\text{uniformly in }B_{1/4}(\nu^0),
\end{align*}
for some $z^0\in\overline{B_3}$ and $\bar\nu^0\in\p B_1$. Note that from $\bar\nu^i_1\le0$, $\bar\nu^0_1\le 0$.

In fact, we have the above convergence in the $C^1$-sense as well. For this purpose, we observe that by the asymptotic property, there exist constants $c_1>0$ and $C_1>0$, independent of $i$, such that $$
c_1\le\frac{u}{\rho_i^\be}\le C_1,\quad c_1\le\frac{u_{\tau_i}}{\bar\rho_i^\be}\le C_1\quad\text{in }B_{\rho_i/2}(x^i).
$$
This is equivalent to $$
c_1\le q_i\le C_1,\quad c_1\le \left(\frac{\rho_i}{\bar\rho_i}\right)^\be \bar q_i\le C_1\quad\text{in }B_{1/2}(\nu^i).
$$
From $1\le\frac{\bar\rho_i}{\rho_i}\le2$, it follows that $q_i$ and $\bar q_i$ are uniformly bounded below and above by positive constants (independent of $i$) in $B_{1/3}(\nu^0)$. In addition, since $F(D^2u)=u^a$ in $\{u>0\}$ near the free boundary, we have $F(D^2q_i)=q_i^a$ and $F(D^2\bar q_i)=\bar q_i^a$ in $B_{1/3}(\nu^0)$. Thus, for any $\al\in(0,1)$, $$
\|q_i\|_{C^{1,\al}(B_{1/4}(\nu^0))}\le C(n,\la,\Lambda,\al)\left(\|q_i\|_{L^\infty(B_{1/3}(\nu^0))}+\|q_i^a\|_{L^\infty(B_{1/3}(\nu^0))}\right)\le C,
$$
and similarly, $$
\|\bar q_i\|_{C^{1,\al}(B_{1/4}(\nu^0))}\le C.
$$
Therefore, \begin{align*}
    &q_i\to q_0\quad\text{in }C^1_{\loc}(B_{1/4}(\nu^0)),\\
    &\bar q_i\to \bar q_0\quad\text{in }C^1_{\loc}(B_{1/4}(\nu^0)).
\end{align*}

Now, to reach a contradiction (and complete the proof), we recall that $\D w_{\tau_i}(x^i)=0$, which is equivalent to $\D(q_i-\bar q_i)(\nu^i)=0$. By the $C^1$-convergence we have $\D q_0(\nu^0)=\D\bar q_0(\nu^0)$. Comparing their first components, we get (by using $q_0(x)=A_0(x\cdot\nu^0)_+^\be$ and $\bar q_0(x)=A_0((x+z^0)\cdot\bar\nu^0)_+^\be$) $$
\nu^0_1(\nu^0\cdot\nu^0)_+^{\be-1}=\bar\nu^0_1((\nu^0+z^0)\cdot\bar\nu^0)^{\be-1}_+.
$$
Since $\nu^0_1>0$ and $\bar\nu_1^0\le0$, we see that the left-hand side of the equation is strictly positive, while the right-hand side is nonpositive. This is a contradiction.


\section{Proof of Theorem~\ref{thm:sym-FN-par}}
In analogy with the elliptic case, we use the following notations in this proof: for $0\le \tau<1$ and $z=(x,t)=(x_1,\cdots,x_n,t)$, \begin{align*}
    \begin{cases}
        \Pi_\tau=\{z\in\R^{n+1}\,:\,x_1=\tau\}\,:\,\text{the hyperplane},\\
        z^\tau=(2\tau-x_1,x_2,x_3,\cdots\,x_n,t)\,:\,\text{the reflection of }z\text{ with respect to }\Pi_\tau,\\
        \Sigma_\tau=\{z\,:\,z^\tau\in Q_1,\,x_1<\tau\}\,:\,\text{the reflection of the right side of $Q_1$},\\
        u_\tau(z)=u(z^\tau)\,:\,\text{the reflection of }u\text{ with respect to }\Pi_\tau.
     \end{cases}   
\end{align*}
We observe that $u_\tau$, $0<\tau<1$, satisfies the equation 
\begin{align}
    \label{eq:u_tau-pde-par}
    F(D^2u_\tau)-\p_tu_\tau=f(u_\tau)\quad\text{in }\Sigma_\tau\cap\{u_\tau>0\},
\end{align}
and the asymptotic property at each $z^0=(x^0,t^0)\in\p\{u_\tau>0\}$ with $t^0>T_1$ \begin{align}
    \label{eq:u_tau-asymp-par}
    \left|u_\tau(x,t)-A_{(z^0)^\tau}\left((x-x^0)\cdot\tilde\mu^{(z^0)^\tau}\right)_+^\be\right|\le C_0\left(|x-x^0|+\sqrt{|t-t^0|}\right)^{2+\de_\be},\quad(x,t)\in\tilde Q_{d_{(z^0)^\tau}}(z^0).
\end{align}
They follow from Hessian property of $F$ and the asymptotic property \eqref{eq:asymp-par-x} of $u$ at $(z^0)^\tau\in\p\{u>0\}$, respectively, as in their elliptic counterparts \eqref{eq:u_tau-pde} and \eqref{eq:u_tau-asymp}.

We first prove the symmetry and monotonicity results in Theorem~\ref{thm:sym-FN-par} for $0<t<T_1$, where there is no free boundary points.
 As in the elliptic case, by Hessian property of $F$, it is enough to show that $u_\tau\ge u$ in $\Sigma_\tau$ for every $0<\tau<1$. 
To prove it by contradiction, we assume that the set $D^{T}_\tau:=\{u_\tau<u\}\cap\Sigma_\tau\cap\{0<t<T\}$ is nonempty for some $0<\tau<1$ and $0<T<T_1$. 

By using the additional assumption on $f$ in Theorem~\ref{thm:sym-FN-par}, \eqref{eq:u_tau-pde-par} and $\inf_{B_1\times[0,T]}u>0$, $w_\tau:=u_\tau-u  $  satisfies in $D^{T}_\tau = \{w_\tau < 0\}$ 

\begin{align*}
    P^-_{\la,\Lambda}(D^2w_\tau)-\p_tw_\tau&\le (F(D^2u_\tau)-\p_tu_\tau)-(F(D^2u)-\p_tu)\\
    &=f(u_\tau)-f(u)=-(f(u)-f(u_\tau))\\
    &\le \ka_0 u_{\tau}^{a-1}(u-u_\tau)\le -\ka_0\left(\inf_{B_1\times[0,T]}u\right)^{a-1}w_\tau.
\end{align*}
Moreover, \eqref{eq:sol-PDE-FN-par} implies that $w_\tau\ge0$ on $\p_p\Sigma_\tau$, thus $w_\tau=0$ on $\p_p D_\tau^{T}$. For $m:=\ka_0\left(\inf_{B_1\times[0,T]}u\right)^{a-1}$, we define $W_\tau(x,t):=e^{-m(t-t^0)}w_\tau(x,t)$. Then
$$
W_\tau=0\quad\text{on }\p_pD^{T}_\tau, \qquad W_\tau < 0  \quad \hbox{in }\ D^{T}_\tau,
$$
 and \begin{align*}
     P^-_{\la,\Lambda}(D^2W_\tau)-\p_tW_\tau=e^{-m(t-t^0)}\left( P^-_{\la,\Lambda}(D^2w_\tau)-\p_tw_\tau\right)+me^{-m(t-t^0)}w_\tau\le 0\quad\text{in }D^{T}_\tau .
\end{align*}
Note that in $D_\tau^T$, $W_\tau$ takes its minimum at an interior point $z^0=(x^0,t^0)$, say. We then take $r>0$ so that the cylinder $Q:=B_r(x^0)\times(t^0-r^2,t^0)$ is inside $D_\tau^T$ and $\p_pQ$ touches $\p_pD^T_\tau$. Then the minimum principle implies that $W_\tau\equiv\inf_{D_\tau^T}W_\tau<0$ in $Q$. This is a contradiction since $W_\tau=0$ on $\p_pD_\tau^T$.

By continuity $u(\cdot,T_1)$ is radial symmetric and $\p_{|x|}u(\cdot,T_1)\ge 0$ in $B_1$, which implies $\{u(\cdot,T_1)=0\}$ is either a one-point set $\{\mathbf0\}$ or a ball $B_{r_0}$ for some $0<r_0<1$. We want to show $\{u(\cdot,T_1)=0\}=\{\mathbf0\}$. For this aim, and towards a contradiction,  assume $\{u(\cdot,T_1)=0\}=B_{r_0}$. From $(r_0e_1,T_1)\in\p\{u>0\}$, we can find a small constant $0<\rho\ll r_0$ such that $u<\e_1$ in a cylinder $B_{3\rho}(r_0e_1)\times(T_1-\rho^2,T_1]$, where $\e_1$ is as in \eqref{eq:assump-f-FN}. We consider a smaller  cylinder $A:=B_\rho((r_0-2\rho)e_1)\times(T_1-\rho^2,T_1]$ and let $\tau:=r_0-\rho\in(0,1)$. Then $w_\tau=u_\tau-u\ge0$ in $A\subset B_1\times(0,T_1]$ and 
$$P^-_{\la,\Lambda}(D^2w_\tau)-\p_tw_\tau\le f(u_\tau)-f(u)=u_\tau^a-u^a\le0 \quad  \hbox{ in } A.$$
 Since $w_\tau((r_0-2\rho)e_1,T_1)=u(r_0e_1,T_1)-u((r_0-2\rho)e_1,T_1)=0-0=0$ and $((r_0-2\rho)e_1,T_1)\in A\setminus\p_pA$, the minimum principle gives $w_\tau\equiv0$ in $A$. However, using $\{u(\cdot,T_1)=0\}=B_{r_0}$, we have $w_\tau((r_0-3\rho)e_1,T_1)=u((r_0+\rho)e_1,T_1)-u((r_0-3\rho)e_1,T_1)>0$ for a point $((r_0-3\rho)e_1,T_1)\in\overline A$. This is a contradiction.

We now prove Theorem~\ref{thm:sym-FN-par} for every $0<t<1$. As before, it suffices to show that $u_\tau\ge u$ in $\Sigma_\tau$ for each $0<\tau<1$. Since $F(D^2u)-\p_tu=f(u)\ge0$ near $\p_pQ_1$, we have by parabolic Hopf Lemma (Theorem~4.2 in \cite{CafLiNir13}) that \begin{align}
    \label{eq:par-Hopf-bdry}
    \p_\nu u<0\quad\text{on }\p B_1\times(0,1],
\end{align}
for any spatial unit vector $\nu$ on $\p B_1\times(0,1]$ which points into $Q_1$. This implies that $\p_{x_1}u>0$ in a neighborhood of $\{e_1\}\times[T_1/2,1]$. Combining this with the monotonicity of $u(\cdot,t)$ for $t<T_1$ gives that for any $\tau\in(0,1)$ close to $1$, we have $u_\tau\ge u$ in $\Sigma_\tau$. Hence we can start our moving plane and we let $\tau_0$ to be the smallest value such that $u_\tau\ge u$ in $\Sigma_\tau$ for all $\tau_0<\tau<1$.

\indent Towards a contradiction, we assume $\tau_0>0$, and take sequences $\e_i\searrow0$ and $\tau_i:=\tau_0-\e_i>0$ such that $D_i:=\{w_{\tau_i}<0\}\cap\Sigma_{\tau_i}\neq\emptyset$, where $w_{\tau_i}:=u_{\tau_i}-u$. For such $D_i$'s, we define $$
D_0:=\{z^0\in\overline{Q_1}\,:\, z^{i_k}\to z^0\text{ for some subsequence }z^{i_k}\in D_{i_k}\}.
$$
Note that $D_i\neq\emptyset$ implies $D_0\neq\emptyset$. For $w_{\tau_0}:=\lim_{i\to\infty}w_{\tau_i}=u_{\tau_0}-u$, we claim that $$
w_{\tau_0}=|\D w_{\tau_0}|=0\quad\text{in }D_0.
$$
Indeed, from $w_{\tau_i}<0$ in $D_i$ we have $w_{\tau_0}\le 0$ on $D_0$. Moreover, by the definition of $\tau_0$, we also have $w_{\tau_0}\ge 0$ in $\overline{\Sigma_{\tau_0}}\supset D_0$, thus $w_{\tau_0}=0$ on $D_0$. To see that $\D w_{\tau_0}=0$ on $D_0$, let $z^0\in D_0$ and $z^i\in D_i$ with $z^i\to z^0$ over a subsequence. Fix a spatial unit vector $e\in\p B_1$. Since each $D_i$ is open relative to $Q_1$, there is a line segment $(y_1^i,y_2^i)\subset D_i$ such that $y_1^i$, $y_2^i\in\p_p D_i$, $y_1^i=y_2^i+se$ for some $s>0$, and $z^i$ lies on $(y^i_1,y^i_2)$. Since $w_{\tau_i}\ge 0$ on $\p_p \Sigma_{\tau_i}$, we have $w_{\tau_i}=0$ on $\p_p D_i$. In particular, $w_{\tau_i}(y_1^i)=w_{\tau_i}(y_2^i)=0$. Since $w_{\tau_i}$ is at least pointwise $C^1_x$ in $\overline{\Sigma_{\tau_i}}$, if $y^i_3\in(y^i_1,y^i_2)$ is a local minimum point of $w_{\tau_i}$ on $(y^i_1,y^i_2)$, then $\p_ew_{\tau_i}(y^i_3)=0$. Then, over a subsequence $[z^i,y_3^i]$ converges either to a point $\{z^0\}$ or to a line segment $[z^0,y^0]\subset \overline{D_0}$. In either case we have $\p_ew_{\tau_0}(z^0)=0$. Since $e$ and $z^0$ are arbitrary, we see that $\D w_{\tau_0}=0$ on $D_0$.

\indent We now decompose $D_0$ into two parts $$
D_0=(D_0\cap\{u>0\})\cup(D_0\cap\{u=0\})
$$
and will prove that these sets are empty, contradicting $D_0\neq\emptyset$.\\


\noindent
\underline{\emph{Claim  A:}}\quad
 $D_0\cap\{u>0\}=\emptyset.$ \\

Note that \emph{Claim A} follows once we show \\
(A1) $D_0\cap\Sigma_{\tau_0}\cap\{u>0\}=\emptyset$,\\
(A2) $D_0\cap \Pi_{\tau_0}\cap\{u>0\}=\emptyset.$\\

Since $w_{\tau_0}=0$ in $D_0$, to prove (A1) it is enough to show that \begin{align}
    \label{eq:w>0-par}
    w_{\tau_0}>0\quad\text{in }\Sigma_{\tau_0}\cap\{u>0\}.
\end{align}
To this aim, we recall that $u_{\tau_0}\ge u$ in $\Sigma_{\tau_0}$ and observe that in $\Sigma_{\tau_0}\cap\{u>0\}$ 
\begin{align}\label{eq:w-subsol-1-par}\begin{split}
P^-_{\la,\Lambda}(D^2w_{\tau_0})-\p_tw_{\tau_0}&\le (F(D^2u_{\tau_0})-\p_tu_{\tau_0})-(F(D^2u)-\p_tu)\\
&=f(u_{\tau_0})-f(u)\le \ka_0u^{a-1}w_{\tau_0}.
\end{split}\end{align}
Let $C$ be a connected component of $\Sigma_{\tau_0}\cap\{u>0\}$. Then, using \eqref{eq:w-subsol-1-par} and the fact that $u^{a-1}$ is bounded in every compact subset of $\Sigma_{\tau_0}\cap\{u>0\}$, we can deduce by parabolic Hopf Lemma that $w_{\tau_0}$ cannot attain its minimum in the interior of $C$ unless it is a constant in $C$. In particular, we have either $w_{\tau_0}\equiv0$ in $C$ or $w_{\tau_0}>0$ in $C$. To prove \eqref{eq:w>0-par} by contradiction, we assume that $w_{\tau_0}\equiv0$ in some component $C$. From the fact that $\Sigma_{\tau_0}\cap\{u>0\}$ is relatively open in $Q_1$ and $w_{\tau_0}>0$ on $(\p_p\Sigma_{\tau_0}\setminus\Pi_{\tau_0})\cap Q_1$, it follows that $C$ is relatively open and $\p_p C\cap(\p_p\Sigma_{\tau_0}\setminus\Pi_{\tau_0})\cap Q_1=\emptyset$. This enables us to take two points $y\in\p_p C\cap\Sigma_{\tau_0}$ and $z\in C$ such that $u(y)=u_{\tau_0}(y)=0$, $u(z)=u_{\tau_0}(z)>0$, and $z-y=(r e_1,0)$ for some $r>0$. Then, the reflected points $y^{\tau_0}$, $z^{\tau_0}\in Q_1\cap\{x_1>\tau_0\}$ satisfy that $u(y^{\tau_0})=0$, $u(z^{\tau_0})>0$, and $y^{\tau_0}-z^{\tau_0}=(r e_1,0)$. This is a contradiction, since $u$ is nondecreasing in $x_1$-direction in $Q_1\cap\{x_1>\tau_0\}$ by the definition of $\tau_0$.\\

To prove (A2) by contradiction, we suppose that there exists a point $z^0=(x^0,t^0)\in D_0\cap\Pi_{\tau_0}\cap\{u>0\}$. Note that $t^0\ge T_1$. If $x^0\in \p B_1$, then $\p_{x_1}u(z^0)=-\frac12\p_{x_1}w_{\tau_0}(z^0)=0$, which contradicts \eqref{eq:par-Hopf-bdry}. Thus we may assume $x^0\in B_1$. Then, we can find a small $r>0$ such that a cylinder $B_r(x^0-re_1)\times(t^0-r^2,t^0)$ is contained in $\Sigma_{\tau_0}\cap\{u>0\}$. Notice that the cylinder touches $\Pi_{\tau_0}$ on $\{x^0\}\times(t^0-r^2,t^0)$. By taking smaller $r$ if necessary, we may assume that $u>c_0$ in the cylinder for some positive constant $c_0>0$. From \eqref{eq:w>0-par} and \eqref{eq:w-subsol-1-par}, together with the fact that $w_{\tau_0}=0$ on $\Pi_{\tau_0}$, we have $\p_{x_1}w_{\tau_0}(z^0)<0$ by Hopf Lemma. This contradicts that $\D w_{\tau_0}=0$ on $D_0$, and completes the the proof of \emph{Claim A}.\\


\noindent
\underline{\emph{Claim  B:}}\quad
 $D_0\cap\Sigma_{\tau_0}\cap\{u=0\}=\emptyset$.\\

For each open set $D_i=\{w_{\tau_i}<0\}\cap\Sigma_{\tau_i}$ defined in the beginning of the proof, since $w_{\tau_i}<0$ in $D_i$ and $w_{\tau_i}\ge 0$ on $\p_pD_{i}$, we can take a point $z^i=(x^i,t^i)\in D_i$ such that $\D w_{\tau_i}(z^i)=0$. Then, over a subsequence $z^i\to z^0=(x^0,t^0)\in D_0$. From \emph{Claim A}, $z^0\in D_0\cap\{u=0\}$. Note that from $w_{\tau_0}=u_{\tau_0}-u=M_t-u>0$ on $(\p_p\Sigma_{\tau_0}\setminus\Pi_{\tau_0})\cap Q_1$, we have $D_0\subset \Sigma_{\tau_0}\cup\Pi_{\tau_0}$. Let $\nu^0=(\nu^0_x,\nu^0_t)=(\nu^0_1,\cdots,\nu^0_n,\nu^0_t)$ be the normal vector to $\p\{u>0\}$ at $z^0$ pointing toward $\{u>0\}$. We then have the following three possibilities:
\begin{enumerate}[(i)]
\item $z^0\in D_0\cap\Sigma_{\tau_0}\cap\{u=0\}$.
\item $z^0\in D_0\cap\Pi_{\tau_0}\cap\{u=0\}$ and $\Pi_{\tau_0}$ is orthogonal to $\p\{u>0\}$ at $z^0$ (i.e., $\nu^0_1=0$, $\nu^0_x\neq0$).
\item $z^0\in D_0\cap\Pi_{\tau_0}\cap\{u=0\}$ and $\Pi_{\tau_0}$ is non-orthogonal to $\p\{u>0\}$ at $z^0$ (i.e., $\nu^0_1\neq0$).
\end{enumerate}

Before discussing the above three cases, we prove that $z^0\in\p\{u>0\}\cap\p\{u_{\tau_0}>0\}$ and $t^0>T_1$.

Indeed, we easily have $z^0\in\p\{u>0\}$ due to $u(z^i)>u_{\tau_i}(z^i)>0$ and $u(z^0)=0$. Moreover, it follows from $u>0$ on $\p_pQ_1$ that $D_0\cap\{u=0\}\subset \Sigma_{\tau_0}\cup(\Pi_{\tau_0}\setminus\p_pQ_1)$, which implies $z^i\in\Sigma_{\tau_0}$ for large $i$. Using $u_{\tau_0}\ge u$ in $\Sigma_{\tau_0}$ we also have $u_{\tau_0}(z^i)\ge u(z^i)>0$. Those combined with $u_{\tau_0}(z^0)=w_{\tau_0}(z^0)+u(z^0)=0$ yield $z^0\in\p\{u_{\tau_0}>0\}$. It remains to prove $t^0>T_1$. Otherwise, we have $t^0=T_1$, since $z^0$ is a free boundary point. Recalling $\{u(\cdot,T_1)=0\}=\{\mathbf0\}$, we see that $z^0=(\mathbf0,T_1)$. From $\tau_0>0$, we obtain $w_{\tau_0}(z^0)=u\left((z^0)^{\tau_0}\right)-u(z^0)>0$, which contradicts $w_{\tau_0}=0$ on $D_0$, as desired.

\medskip

We deal with the cases (i) and (ii) at the same time, as they both follow from the application of the result in Appendix~\ref{appen:comp-prin}. As in the elliptic problem, the (parabolic) asymptotic properties \eqref{eq:asymp-par-x} and \eqref{eq:u_tau-asymp-par} give  \begin{align}
    \label{eq:w-small-par}
    w_{\tau_0}(x,t)=O\left(\left(|x-x^0|+\sqrt{|t-t^0|}\right)^{2+\de_\be}\right).
\end{align}
For small $r>0$, $$
P^-_{\la,\Lambda}(D^2w_{\tau_0})-\partial_tw_{\tau_0}\le f(u_{\tau_0})-f(u)=u_{\tau_0}^a-u^a\le0\quad\text{in }A_r:=\Sigma_{\tau_0}\cap\{u>0\}\cap\tilde Q_r(z^0).
$$
Let $(n+1)$-dimensional domain $\Omega^\e$ and the function $H^\e$ be as in Appendix \ref{subsec:comparison-par}, with small $\e>0$ satisfying \eqref{eq:comparison-L-par}. Due to the $C^{1,|a|}_x\cap C^{1/2+|a|/2}_t$-regularity of $\p\{u>0\}$, we can take possible rotations and translations on $\Omega^\e$ and $H^\e$ to get a domain $\Omega$ and a function $H$ as well as a sequence $z^j\in \Omega\cap\tilde Q_r(z^0)$ with $z^j\to z^0$ such that for small $r>0$ \begin{align}\label{eq:H-Omega}\begin{cases}
    \Omega\cap\tilde Q_r(z^0)\subset A_r,\\
    P_{\la,\Lambda}^-(D^2H)-\p_tH=0\quad\text{in }\Omega\cap\tilde Q_r(z^0),\quad H=0\quad\text{on }\p_p\Omega\cap\tilde Q_r(z^0),\\
    H(z^j)\ge|x^j_1-x^0_1|^{2+\sigma},\quad w_{\tau_0}(x^j)\le C|x^j_1-x^0_1|^{2+\de_\be}.
\end{cases}\end{align}
Here, the last inequality containing $w_{\tau_0}$ follows by using \eqref{eq:w-small-par}. Since $H\le \|h\|_{L^\infty(B_1)}$ on $\p_p\Omega$ and $w_{\tau_0}>0$ on $\overline{\p_p\Omega\setminus\{H=0\}}$, we have for some constant $C^*>0$ $$
H\le C^*w_{\tau_0}\quad\text{on }\p_p\Omega,
$$
and obtain, by applying the comparison principle, $$
H\le C^*w_{\tau_0}\quad\text{in }\Omega.
$$
By taking $\si\in(0,\de_\be)$, this contradicts \eqref{eq:H-Omega}.

\medskip 

We now consider the case (iii). The assumption $\nu^0_1\neq0$ in (iii) and the fact that $u_{\tau_0}\ge0$ in $\Sigma_{\tau_0}$ imply $\nu^0_1>0$.

Since the normal vector $\nu^0$ at $z^0$ is not parallel to the time axis $(0,\cdots,0,1)$, $\p\{u>0\}\cap\{t=t^0\}$ and $\p\{u>0\}\cap \{t=t^i\}$ are $(n-1)$-dimensional surfaces on $B_1\times\{t^0\}$ and $B_1\times\{t^i\}$, respectively, for large $i$. For such $i$ we take a point $y^i\in\R^n$ so that $z^i_f:=(y^i,t^i)\in\p\{u>0\}\cap\{t=t^i\}$ is closest to $z^i=(x^i,t^i)$. In analogy to notations $\nu^0$ and $\mu^0$, we let $\nu^i=(\nu^i_x,\nu^i_t)=(\nu^i_1,\cdots,\nu^i_n,\nu^i_t)$ be the unit normal to $\p\{u>0\}$ at $z^i_f$ pointing toward $\{u>0\}$, and denote $\mu^i:=\frac{\nu^i_x}{|\nu^i_x|}$. Then, $\mu^i=\frac{x^i-y^i}{\rho_i}$ for $\rho_i:=|x^i-y^i|$. Due to the $C^1$-smoothness of $\p\{u>0\}$ near $z^0$, $\nu^i\to\nu^0$ and $\mu^i\to\mu^0$. In particular, $\nu^i_1\ge\nu^0_1/2>0$ for large $i$. Next, we take a point $\bar y^i\in\R^n$ so that $\bar z^i_f:=(\bar y^i,t^i)\in\p\{u_{\tau_i}>0\}\cap\{t=t^i\}$ is closest to $z^i=(x^i,t^i)$. As before, let $\bar\nu^i=(\bar\nu^i_x,\bar\nu^i_t)=(\bar\nu^i_1,\cdots,\bar\nu^i_n,\bar\nu^i_t)$ be the unit normal to $\p\{u_{\tau_i}>0\}$ at $\bar z^i_f$ pointing toward $\{u_{\tau_i}>0\}$, and denote $\bar\mu^i:=\frac{\bar\nu^i_x}{|\bar\nu^i_x|}$, so that $\bar\mu^i=\frac{x^i-\bar y^i}{\bar\rho_i}$ for $\bar\rho_i:=|x^i-\bar y^i|$. Since $\bar z^i_f\to z^0\in\p\{u_{\tau_0}>0\}$, we have $\bar\nu^i\to\bar\nu^0$, where $\bar\nu^0:=(\bar\nu^0_x,\nu^0_t)=(-\nu^0_1,\nu^0_2,\cdots,\nu^0_n,\nu^0_t)$ and $\bar\mu^i\to\bar\mu^0:=\frac{\bar\nu^0_x}{|\bar\nu^0_x|}$. That is, $\bar\mu^0$ is the reflection of $\mu^0$ with respect to $\{x_1=0\}$. Notice that $\rho_i\le\bar\rho_i$, $\bar\nu^i\le0$, and that $z^i$, $z^i_f$, $\bar z_f^i$ $(\bar z_f^i)^{\tau_i}$ all converge to $z^0$.

We claim that $\frac{\bar\rho_i}{\rho_i}\le 2$ for large $i$. Otherwise, we have that over a subsequence $\frac{\bar\rho_i}{\rho_i}>2$. Applying the asymptotic property of $u$ at $z^i_f$ and $u_{\tau_i}$ at $\bar z^i_f$, respectively, we obtain for large $i$ \begin{align*}
    u(z^i)&\le A_{z^i_f}((x^i-y^i)\cdot\mu^i)_+^\be+C_0|x^i-y^i|^{2+\de_\be}=A_{z^i_f}\rho_i^\be+C_0\rho_i^{2+\de_\be}\le \sqrt2A_{z^0}\rho_i^\be,\\
    u_{\tau_i}(z^i)&\ge A_{(\bar z_f^i)^{\tau_i}}((x^i-\bar y^i)\cdot\bar\mu^i)^\be_+-C_0|x^i-\bar y^i|^{2+\de_\be}=A_{\bar z_f^i}\bar\rho_i^\be-C_0\bar\rho_i^{2+\de_\be}\ge\frac{A_{z^0}}{\sqrt2}\bar\rho_i^\be.
\end{align*}
Then, we obtain $$
w_{\tau_i}(z^i)=u_{\tau_i}(z^i)-u(z^i)\ge\sqrt2A_{z^0}\rho_i^\be\left(\frac12\left(\frac{\bar\rho_i}{\rho_i}\right)^\be-1\right)>0,
$$
which contradicts that $z^i\in D_i=\{w_{\tau_i}<0\}\cap\Sigma_{\tau_i}$.

Next, we consider the homogeneous rescalings $$
q_i(x,t):=\frac{u(\rho_ix+y^i,\rho_i^2t+t^i)}{\rho_i^\be},\qquad \bar q_i(x,t):=\frac{u_{\tau_i}(\rho_ix+y^i,\rho_i^2t+t^i)}{\rho_i^\be}.
$$
Again, we use the asymptotic property of $u$ at $z^i_f$ and $u_{\tau_i}$ at $\bar z^i_f$, respectively, to get \begin{align*}
    &|q_i(x,t)-A_{z_f^i}(x\cdot\mu^i)^\be_+|\le C_0\rho_i^{2+\be_\be-\be}\left(|x|+|t|^{1/2}\right)^{2+\de_\be},\\
    &\left|\bar q_i(x,t)-A_{(\bar z_f^i)^{\tau_i}}\left(\left(x+\frac{y^i-\bar y^i}{\rho_i}\right)\cdot\bar\mu^i\right)^\be_+\right|\le C_0\rho_i^{2+\de_\be-\be}\left(\left|x+\frac{y^i-\bar y^i}{\rho_i}\right|+|t|^{1/2}\right)^{2+\de_\be}.
\end{align*}
For large $i$, these estimates hold in $\tilde Q_{1/2}(\mu^i,0)$, and thus in $\tilde Q_{1/4}(\mu^0,0)$ as well. From $\frac{\bar\rho_i}{\rho_i}\le2$, we see that $$
\frac{|y^i-\bar y^i|}{\rho_i}\le\frac{|y^i-x^i|+|x^i-\bar y^i|}{\rho_i}=\frac{\rho_i+\bar\rho_i}{\rho_i}\le 3.
$$
Thus, over a subsequence, \begin{align*}
    &q_i\to q_0(x,t):=A_{z^0}(x\cdot\mu^0)_+^\be\quad\text{uniformly in }B_{1/4}(\mu^0,0),\\
    &\bar q_i\to \bar q_0(x,t):=A_{z^0}((x+\eta^0)\cdot\bar\mu^0)_+^\be\quad\text{uniformly in }B_{1/4}(\mu^0,0)
\end{align*}
for some $\eta^0\in \overline{B_3}$.

In fact, we have the above convergence in the $C^1$-sense as well. For this purpose, we observe that by the asymptotic property, there exist constants $C_1>c_1>0$, independent of $i$, such that $$
c_1\le \frac{u}{\rho_i^\be}\le C_1,\quad c_1\le\frac{u_{\tau_i}}{\bar\rho_i^\be}\le C_1\quad\text{in }\tilde Q_{\rho_i/2}(z^i).
$$
This is equivalent to $$
c_1\le q_i\le C_1,\quad c_1\le\left(\frac{\rho_i}{\bar\rho_i}\right)^\be\bar q_i\le C_1\quad\text{in }\tilde Q_{1/2}(\mu^i,0).
$$
From $1\le\frac{\bar\rho_i}{\rho_i}\le2$, it follows that $q_i$ and $\bar q_i$ are uniformly bounded above and below by positive constants, independent of $i$, in $\tilde Q_{1/3}(\mu^0,0)$. In addition, since $F(D^2u)-\p_tu=u^a$ in $\{u>0\}$ near the free boundary, we have $F(D^2q_i)-\p_tq_i=q_i^a$ and $F(D^2\bar q_i)-\p_t\bar q_i=\bar q_i^a$ in $\tilde Q_{1/3}(\mu^0,0)$. Thus, for any $\al\in(0,1)$,\begin{align*}
    \|q_i\|_{C^{1,\al}_x(\tilde Q_{1/4}(\mu^0,0))}\le C(n,\la,\Lambda,\al)\left(\|q_i\|_{L^\infty(\tilde Q_{1/3}(\mu^0,0))}+\|q_i^a\|_{L^\infty(\tilde Q_{1/3}(\mu^0,0))}\right)\le C,
\end{align*}
and similarly, $$
\|\bar q_i\|_{C^{1,\al}_x(\tilde Q_{1/4}(\mu^0,0))}\le C.
$$
Therefore, \begin{align*}
    &q_i\to q_0\quad\text{in }C^1_{\loc}(\tilde Q_{1/4}(\mu^0,0)),\\
    &\bar q_i\to\bar q_0\quad\text{in }C^1_{\loc}(\tilde Q_{1/4}(\mu^0,0)).
\end{align*}

To reach a contradiction, we recall that $\D w_{\tau_i}(z^i)=0$, which implies $\D(q_i-\bar q_i)(\mu^i,0)=0$. By the $C^1$-convergence we have $\D q_0(\mu^0,0)=\D\bar q_0(\mu^0,0)$. Comparing their first components, we get $$
\mu^0_1(\mu^0\cdot\mu^0)^{\be-1}_+=-\mu^0_1((\mu^0+\eta^0)\cdot\bar\mu^0)^{\be-1}_+,
$$
where $\mu^0_1$ is the first component of $\mu^0$. From $\nu^0_1>0$, we have $\mu^0_1>0$, thus the left-hand side of the equation is strictly positive, while the right-hand side is nonpositive. This is a contradiction.



\appendix

\section{Properties of $A_{x^0}$ and $A_{z^0}$}\label{sec:appen-A}

In this section we discuss values and properties of $A_{x^0}>0$ and $A_{z^0}>0$, when $q_{x^0}(x)=A_{x^0}(x\cdot\nu^{x^0})_+^\be$ and $q_{z^0}(x,t)=A_{z^0}(x\cdot\mu^{z^0})_+^\be$ are solutions to \begin{align*}
    &F(D^2q_{x^0})=q_{x^0}^a\quad\text{in }B_1,\\
    &F(D^2q_{z^0})-\p_tq_{z^0}=q_{z^0}^a\quad\text{in }\tilde Q_1.
\end{align*}
Here, $\nu^{x^0}$ and $\mu^{z^0}$ are unit vectors in $\R^n$. To get $A_{x^0}$, we compute \begin{align*}
    &F(D^2q_{x^0})=A_{x^0}\be(\be-1)(x\cdot\nu^{x^0})_+^{\be-2}F(\nu_{x^0}\otimes\nu_{x^0}),\\
    &q_{x^0}^a=A_{x^0}^a(x\cdot\nu^{x^0})_+^{\be a}=A_{x^0}^a(x\cdot\nu^{x^0})_+^{\be-2}.
\end{align*}
They readily give $$
A_{x^0}=[\be(\be-1)F(\nu^{x^0}\otimes\nu^{x^0})]^{-\be/2}.
$$
For its bounds, we observe that $\nu^{x^0}\otimes\nu^{x^0}$ has the eigenvalues $1,0,\cdots,0$, which implies $\la\le F(\nu^{x^0}\otimes\nu^{x^0})\le\Lambda$. Thus, we obtain the uniform lower and upper bounds on $A_{x^0}$ $$
[\be(\be-1)\Lambda]^{-\be/2}\le A_{x^0}\le [\be(\be-1)\la]^{-\be/2}.
$$
Since $q_{z^0}$ is time independent, it also satisfies $F(D^2q_{z^0})=q_{z^0}^a$. Thus, repeating the above process with $q_{z^0}$ will yield $$
A_{z^0}=[\be(\be-1)F(\mu^{z^0}\otimes\mu^{z^0})]^{-\be/2}
$$
and $$
[\be(\be-1)\Lambda]^{-\be/2}\le A_{z^0}\le [\be(\be-1)\la]^{-\be/2}.
$$
Finally, if $\nu$, $\mu\in\p B_1$ with $F(\nu\otimes\nu)\ge F(\mu\otimes\mu)$, then for the eigenvalues $e_i$'s of $\nu\otimes\nu-\mu\otimes\mu$, \begin{align*}
    F(\nu\otimes\nu)-F(\mu\otimes\mu)&\le P^+_{\la,\Lambda}(\nu\otimes\nu-\mu\otimes\mu)\le \Lambda\Sigma|e_i|\\
    &\le n\Lambda\|\nu\otimes\nu-\mu\otimes\mu\|\le C(n,\Lambda)\|\nu-\mu\|.
\end{align*}
Therefore, $\nu\mapsto A_\nu$ is continuous.


\section{A non-standard comparison principle}\label{appen:comp-prin}

For small $\si>0$, we construct\footnote{After finalizing this paper we found out that a different (but similar  in nature) construction for elliptic case has been done by Silvestre-Sirakov \cite{SilSir15}, that hinges on an earlier result by Armstrong-Sirakov-Smart \cite{ASS2012}, where an  higher dimensional homogeneous solution in cones are constructed. }
 in Appendix~\ref{subsec:comparison-ell} an $n$-dimensional Lipschitz cone $K$, centered at the origin  and contained in $\{x_1>0,\, x_2>0\}$, and a function $H:K\to\R$ satisfying the following: for $1\le \Lambda/\la\le \eta_0$ with $\eta_0>1$ depending only on $\si$, \begin{align*}
    \begin{cases}
        \text{- }P^-_{\la,\Lambda}(D^2H)=0,\quad H\ge0\,\,\,\text{in }K,\quad H=0\,\,\,\text{on }\p K,\\
        \text{- }H\text{ does not decay faster than }|x|^{2+\si}\text{ near the center of }K,\text{ see }\eqref{eq:comparison-L}.
    \end{cases}
\end{align*}
We prove its parabolic counterpart in Appendix~\ref{subsec:comparison-par}, see \eqref{eq:comparison-L-par}.


\subsection{Elliptic Case}\label{subsec:comparison-ell}
We fix a small constant $\sigma_0\in(0,\si)$, and consider a function $\tilde h$, defined in a $2$-dimensional cone $\mathcal C$,  such that
\begin{align}\label{eq:2dim-hom-sol-cone}
\begin{cases}
    \text{- }\mathcal C\subset\{x_1>0,\,x_2>0\}\subset\R^2\text{ with opening }\pi/2-\mu,\\
    \text{- } \tilde h\,:\, \mathcal C\to\R\text{ is a positive solution to }P^-_{\la,\Lambda}(D^2\tilde h)=0\text{ in }\mathcal C\text{ and }\tilde h=0\text{ on }\p \mathcal C,\\
    \text{- }\tilde h\text{ is homogeneous of degree }2+\si_0.
\end{cases}
\end{align}
For $1\le\Lambda/\la\le \eta_0=\eta_0(\si_0)$ and small $\mu>0$, this follows from \cite{Leo17}. Indeed, thanks to Theorem~2.4 in \cite{Leo17} we can find a nonnegative $\al$-homogeneous solution $\tilde h$ of $P^-_{\la,\Lambda}(D^2\tilde h)=0$ in $\mathcal C$ with zero boundary value. Here, the homogeneity $\al$ is determined from the equation $$
g_\eta(\al)=\frac12(\pi/2-\mu),
$$
where $$g_\eta(\al)=\arctan\sqrt\eta+\frac{2-\al}{\sqrt{(\al-1+1/\eta)(\al-1+\eta)}}\arctan\sqrt{\frac{\eta(\al-1+1/\eta)}{\al-1+\eta}},\quad \eta=\Lambda/\la.  $$
We can see that $\al$ is completely determined by the opening of the cone, $\pi/2-\mu$, and the value $\eta=\Lambda/\la$. To see if we can make $\al=2+\si_0$, we first take $\eta=1$ and compute $g_1(\al)=\frac\pi4\left(\frac1{\al-1}\right)$. Solving $\frac\pi4\left(\frac1{\al-1}\right)=\frac12(\pi/2-\mu)$ gives $\al=2+\frac{2\mu}{\pi-2\mu}$, thus $\al=2+\si_0$ holds for small $\mu>0$ (specifically, $\mu=\frac{\si_0\pi}{2(\si_0+1)}$. This computation implies that for $\eta=\Lambda/\la$ larger than $1$, say $1<\eta<\eta_0=\eta_0(\si_0)$, we can still find small $\mu>0$ to have $\al=2+\si_0$.

\medskip

Next, by possibly multiplying a constant on $\tilde h$, we assume $|x|^{2+\si}< \tilde h/2$ on the line $\{0<x_1=x_2<1\}$. We then extend $\tilde h$ from $\mathcal C$ to the cylindrical cone $\mathcal C\times\R^{n-2}$  by defining 
 $$
h(x_1,x_2,x''):= \tilde h(x_1,x_2),\quad(x_1,x_2)\in\mathcal C,\,\, x''\in\R^{n-2}.
$$
For small $\e>0$ we define $n$-dimensional cones
$$
K^\e:=\{h>0\}\cap\{x_1+x_2>\e|x''|\}, \qquad K^0 = \{h > 0\},
$$
and observe that $\lim_{\e \to 0}K^\e = K^0$.

Let now $H^\e$ be a solution of \begin{align*}
\begin{cases}
    \text{- }P^-_{\la,\Lambda}(D^2H^\e)=0\text{ in }K^\e\cap B_1,\\
    \text{- }H^\e=h\text{ on }K^\e\cap\p B_1,\\
    \text{- }H^\e=0\text{ on }\p K^\e\cap B_1.
\end{cases}
\end{align*}
By applying the comparison principle,  $H^\e\le h$  in $K^\e\cap B_1$, and moreover 
\begin{equation}\label{eq:convergence}
H^\e \to h\quad\text{in } C^1_{\loc}(K^0\cap B_1), \quad \hbox{ and } \quad H^\e \nearrow  h
\quad \hbox{  pointwise in }K^0\cap B_1.
\end{equation}
We claim now that for $\e >0 $ small enough there is a sequence 
$x^j \in L_{r_j} : = \{x_1=x_2>0,\,x''=0 \} \cap \overline{B_{r_j}} $,  where $r_j := |x^j| \searrow 0$ and such that 
\begin{equation}\label{eq:comparison-L}
H^\e (x^j) \geq  |x^j|^{2+\sigma}.
\end{equation}

To prove  claim \eqref{eq:comparison-L} we argue by contradiction. If \eqref{eq:comparison-L} is not true, then using \eqref{eq:convergence} we can find sequences $\e_j\searrow0$ and $x^j\in L_{r_j}$ with $r_j=|x^j|\to0$ such that \begin{equation}\label{eq:contra}
H^{\e_j}(x^j) = |x^j|^{2+\sigma}, \quad \hbox{and} \quad 
H^{\e_j}(x) <  |x|^{2+\sigma} \quad \forall x\in L_{r_j}.
\end{equation}
Consider 
 $$
\tilde H_{j}(x):=\frac{H^{\e_j}(r_j x)}{(r_j)^{2+\si}},
$$
which satisfies \begin{align*}
    \begin{cases}
        \text{- }P^-_{\la,\Lambda}(D^2\tilde H_{j})=0\quad\text{in }K^{\e_j}\cap B_{1/r_j},\\
        \text{- }\tilde H_{j}=0\quad\text{on }\p K^{\e_j}\cap B_{1/r_j}.
    \end{cases}
\end{align*}
Set  $y^0:=\frac{x^j}{r_j }\in L_1\cap \p B_1$ (note that $L_1\cap \p B_1$ is a singleton). Then 
by construction,  we have 
\begin{equation}\label{eq:control-1}
\tilde H_{j}(y^0) = |y^0|^{2+\si}=1, \qquad  \tilde H_{j}(x ) \leq  |x|^{2+\si} \quad \hbox{for } x\in L_1.
\end{equation}

This, combined with Harnack inequality, implies that in every compact subset of $K^{0}$, \ $\tilde H_j $ is uniformly bounded in $j$. Thus, $\tilde H_{j}\to\tilde H_{0}$ in $K^0$, where $\tilde H_0 $ solves the same PDE as above, and satisfies
\begin{equation}\label{eq:control-2}
\tilde H_{0}(y^0) = 1, \qquad  \tilde H_{0}(x ) \leq  |x|^{2+\si} \quad \hbox{for } x\in L_1.
\end{equation}
Applying  Boundary Harnack Principle to $\tilde H_0$ and $h$ in $K^0 \cap B_1$ implies 
 $$
 c h   \leq  \tilde H_0 
\quad\text{in }K^0\cap B_{1/2}.
 $$
 This in combination with \eqref{eq:control-2} implies 
 $$
 c|x|^{2+ \sigma_0 }\approx  c h   \leq  \tilde H_0  \leq |x|^{2+\sigma}, \qquad x\in L_{1/2} ,
$$
 which is a contradiction, since $\sigma > \sigma_0$.


\subsection{Parabolic case}\label{subsec:comparison-par}
Let the function $h:\R^n\to\R$ and $n$-dimensional cone $K^\e=\{h>0\}\cap \{x_1+x_2>\e|x''|\}$ be as in the elliptic case. For a spatial unit vector $\nu^0_x\in\p B_1$, we define $(n+1)$-dimensional domain $$
\Omega^\e:=\{(x,t)\in\R^{n+1}\,:\,-1<t<1,\,x\in(K^\e\cap B_1)+\e|t|^{1/2}\nu^0_x\}.
$$
Note that every time-slice of $\Omega^\e$ is a truncated cone of shape $K^\e\cap B_1$ with vertex $\e|t|^{1/2}\nu^0_x$. We can observe that as $\e\searrow 0$, $K^\e$ and $\Omega^\e$ converge to $n$-dimensional cylinder $K^0=\{h>0\}$ and $(n+1)$-dimensional cylinder $\Omega^0:=(K^0\cap B_1)\times(-1,1)$, respectively.

Let $H^\e$ be a solution of \begin{align}\label{eq:H-sol}
    \begin{cases}
        \text{- }P^-_{\la,\Lambda}(D^2H^\e)-\p_tH^\e=0\quad\text{in }\Omega^\e,\\
        \text{- For }-1<t<1,\,\,\, H^\e(\cdot,t)=h(\cdot-\e|t|^{1/2}\nu^0_x)\quad\text{on }\p(K^\e\cap B_1)+\e|t|^{1/2}\nu^0_x,\\
        \text{- For }t=-1,\,\,\, H^\e(\cdot,-1)=h(\cdot-\e\nu^0_x)\quad\text{on }(K^\e\cap B_1)+\e\nu^0_x.
    \end{cases}
\end{align}
The Maximum Principle yields $H^\e\le \|h\|_{L^\infty(B_1)}$ in $\Omega^\e$, thus for some function $H^0$ we have over a subsequence 
 \begin{align}
    \label{eq:convergence-par}
    H^\e\to H^0\,\,\,\text{in }C^1_{\loc}(\Omega^0),\,\,\,\text{ and }\,\,H^\e\to H^0\,\text{ pointwise in }\Omega^0.
\end{align}
The fact that $\Omega^0$ is a cylinder gives
 \begin{align*}
    \begin{cases}
        P^-_{\la,\Lambda}(D^2H^0)-\p_tH^0=0\quad\text{in }\Omega^0,\\
        H^0=h\quad\text{on }\p_p\Omega^0,
    \end{cases}
\end{align*}
thus we have by the uniqueness $H^0=h$ in $\Omega^0$.

We now claim that for $\e>0$ small enough there exists a sequence $z^j=(x^j,t^j)=(x^j_1,\cdots,x^j_n,t^j)\in A_{r_j}:=\{(a,a,0,\cdots,0,a^2)\,:\,0<a<r_j\}$, where $r_j\searrow 0$, and such that \begin{align}
    \label{eq:comparison-L-par}
    H^\e(z^j)\ge (x^j_1)^{2+\si}.
\end{align}
To prove \eqref{eq:comparison-L-par} we argue by contradiction. If \eqref{eq:comparison-L-par} is not true, then, using \eqref{eq:convergence-par} we can find sequences $\e_j\searrow0$ and $z^j\in A_{r_j}$ with $r_j=|x^j_1|\to0$ such that $$
H^{\e_j}(z^j)=(x^j_1)^{2+\si},\quad\text{and }\,H^{\e_j}(z)<(x_1)^{2+\si}\text{ for all }z\in A_{r_j}.
$$
Consider $$
\tilde H_j(x,t):=\frac{H^{\e_j}(r_jx,r_j^2t)}{r_j^{2+\si}},
$$
which satisfies \begin{align*}
    \begin{cases}
        P^-_{\la,\Lambda}(D^2\tilde H_j)-\p_p\tilde H_j=0\quad\text{in }\tilde\Omega^j,\\
        \text{For }|t|<1/r_j^2,\,\,\, \tilde H_{j}(\cdot,t)=0\quad\text{on }(\p K^{\e_j}\cap B_{1/r_j})+\e_j|t|^{1/2}\nu^0_x,
    \end{cases}
\end{align*}
where $$
\tilde\Omega^j:=\left\{(x,t)\in\R^{n+1}\,:\, -\frac1{r_j^2}<t<\frac1{r_j^2},\,\,x\in(K^{\e_j}\cap B_{1/r_j})+\e_j|t|^{1/2}\nu^0_x\right\}.
$$
Set $y^0:=(1,1,0,\cdots,0,1)\in \R^{n+1}$. By construction $$
\tilde H_j(y^0)=1,\quad \tilde H_j(z)<(x_1)^{2+\si}\text{ for }z\in A_1.
$$
By parabolic Harnack inequality, in every compact subset of $(K^0\cap B_{1/2})\times(-1/2,1/2)$, $\tilde H^j$ is uniformly bounded in $j$. Thus, over a subsequence $\tilde H^j\to \tilde H^0$ in $(K^0\cap B_{1/2})\times(-1/2,1/2)$, where $\tilde H^0$ satisfies \begin{align*}
    \begin{cases}
        P^-_{\la,\Lambda}(D^2\tilde H^0)-\p_t\tilde H^0=0\quad\text{in }(K^0\cap B_{1/2})\times (-1/2,1/2),\\
        \tilde H^0=0\quad\text{on }(\p K^0\cap B_{1/2})\times(-1/2,1/2),\\
        \tilde H^0(z)\le (x_1)^{2+\si}\text{ for }z\in A_{1/4}.
    \end{cases}
\end{align*}
Applying parabolic Boundary Harnack Principle to $\tilde H^0$ and $h$ in $(K^0\cap B_{1/2})\times(-1/2,1/2)$ implies that $$
ch\le \tilde H^0\quad\text{in }(K^0\cap B_{1/4})\times(-1/4,1/4).
$$
Therefore $$
c(x_1)^{2+\si_0}\approx ch(z)\le \tilde H^0(z)\le (x_1)^{2+\si},\quad z\in A_{1/8},
$$
which is a contradiction since $\si>\si_0$.

\section*{Declarations}

\noindent {\bf  Data availability statement:} All data needed are contained in the manuscript.

\medskip
\noindent {\bf  Funding and/or Conflicts of interests/Competing interests:} The authors declare that there are no financial, competing or conflict of interests.


\end{document}